\documentclass[12pt,a4paper]{article}
\usepackage{graphicx}
\usepackage{graphics}
\usepackage{amsmath, amssymb}
\usepackage[bf]{caption2}
\newcounter{theorem}

\newcounter{theoremcounter}

\newcounter{corollarycounter}
\newtheorem{theorem}[theoremcounter]{Theorem}

\newtheorem{corollary}[corollarycounter]{Corollary}
\usepackage[normalem]{ulem}
\usepackage[cp1251]{inputenc}

\begin{document}

%%%%% TITLE %%%%%

\title{Ergodicity and uniform in time truncation bounds for inhomogeneous birth and death processes with additional transitions from and to origin} %in capital

\date{}
%%%%% AUTHORS %%%%%

\author{
Alexander Zeifman\footnote{Vologda State University, Institute of
Informatics Problems of the  FRC CSC RAS, ISEDT RAS, corresponding
author, telephone/fax +78172721632, e-mail a\_zeifman@mail.ru},
 Anna Korotysheva, Yacov Satin\footnote{Vologda State University, Institute of
Informatics Problems of the  FRC CSC RAS},\\ Rostislav
Razumchik\footnote{Institute of Informatics Problems of the FRC CSC
RAS, Peoples' Friendship University}, Victor Korolev\footnote{
Moscow State University, Institute of Informatics Problems of the
FRC CSC RAS, Moscow, Russia}, Sergey Shorgin\footnote{Institute of
Informatics Problems of the FRC CSC RAS, Moscow, Russia}}

\maketitle

{\bf Abstract.}
In this paper one presents the extension of the
transient analysis of the class
of continuous-time
birth and death processes defined on non-negative integers
with special transitions from and to the
origin. From the origin transitions can occur
to any state. But being in any other state,
besides ordinary transitions to neighbouring states,
a transition to the origin can occur.
All possible transition intensities are assumed to
be non-random functions of time and may depend on the
state of the process. We improve previously known
ergodicity and truncation bounds for this class of processes
which were known only for the case when transitions from the origin decay
exponentially (other intensities must have unique uniform upper bound).
We show how the bounds can be obtained the decay rate is slower than exponential.
Numerical results are also provided.


\section{Introduction}

In this paper consideration is given to one
subclass of continuous-time Markov chains
--
inhomogeneous birth and death processes with additional
transitions from and to origin.
More strictly speaking, one considers
an inhomogeneous
continuous-time Markov chain $\{X(t), \ t\geq 0 \}$
with state space $\mathcal{X}=\{ 0, 1, 2 \dots \}$.
%Denote by $p_{ij}(s,t)=\Pr \left\{ X(t)=j\left| X(s)=i\right. \right\}$,
%$i,j \ge 0, \;0\leq s\leq t$ the transition probabilities of
%$X(t)$.
All possible transition intensities are assumed
to be non-random functions of time and
may depend on the state of the process.
From state $0$ the chain can jump to any state $i>0$
with transition intensity $r_i(t)$.
But transitions from each state
$i>0$ happen with intensities $\mu_{i}(t)$, $\lambda_{i}(t)$ and
$\beta_i(t)$ and can be either to state $(i-1)$ or $(i+1)$ or
$0$, respectively.

Such subclass of processes finds its application in the study of
queueing systems with catastrophes and bulk arrivals (see, for
example, \cite{pk1991,cr1997,cpz2010,cr2004, lc2013,r1,zh,a1,DK,
dzp}). For more details, concerning possible applications, one can
refer to  \cite{a1} and references therein. Note that in the cited
papers authors call transitions from and to origin (governed by
intensities $r_i(t)$ and $\beta_i(t)$) as ``mass arrivals when
empty'' and ``mass exodus'' transitions, respectively. In order to
keep the connection with the previously obtained results in what
follows we use these term as well.

The motivation for this research was given by papers \cite{a1}--\cite{a4},
where authors studied transient behaviour of
various inhomogeneous birth and death processes being a subclass
of $X(t)$.
Specifically their
results concerned ergodicity and perturbation bounds,
and bounds of truncation\footnote{Which allow the calculation of the limiting characteristics
with given accuracy.} of the processes
(for example, probability of being in a particular state at time $t$, or expected
value of the process at time $t$, which started from any given state).
It is well-known that exact computation of transient state probability distribution is not very appealing
way to analyse the behaviour of system described by processes with time-dependent rates.
A good alternative is to have bounds for performance characteristics
of interest, which can be computed fast and are tight enough to make results
meaningful. With respect to this observation, the
direction of research, indicated by papers \cite{a1}--\cite{a4},
looks promising.

Clearly ergodicity, perturbation and truncation bounds heavily
depend on the restrictions imposed on the transition intensities
of the process. The harder restrictions are, the better bounds one can
obtain. But hard restrictions make the results inapplicable
in most cases which are interesting for applications. Thus the main effort
is usually made toward the loosening the restrictions and improving bounds.
In \cite{a1} the authors considered the special case of the process $X(t)$
in which transitions to state 0 cannot depend
on the current state of the process (i.e. $\beta_i(t)=\beta(t)$ for each $i$),
and obtained the first ergodicity, perturbation and truncation bounds.
In this paper by applying a different method
we demonstrate that truncation bounds obtained in papers \cite{a1}
 can be improved. Specifically, it is shown that
one can obtain uniform in time error bounds of truncation
for the general process $X(t)$.
Uniform error bounds of truncation have already been
obtained for several processes in \cite{a2} and \cite{a3}.
In \cite{a3} consideration was given to the process $X(t)$
with $r_i(t)=\beta_i(t)\equiv 0$ (i.e. pure inhomogeneous birth and death processes).
In \cite{a2} one managed to obtain uniform truncation bounds for
the process $X(t)$, in the case when $\beta_i(t)\equiv 0$ and intensities $r_i(t)$ decrease exponentially
with $i$ i.e.
there exists $q>1$ such that $r_i(t) \le q^{-i}$ for any $i>0$.
It is worth noticing that using results from \cite{a2} and \cite{a3},
one can construct truncation bounds for $X(t)$ with minor restrictions
on intensities $r_i(t)$, $\mu_{i}(t)$, $\lambda_{i}(t)$ and
$\beta_i(t)$, but these bounds are uninformative because they
 tend to infinity as  $t \rightarrow \infty$.

In this paper it is shown that one can find
uniform in time error bounds of truncation
for the general process $X(t)$ under less
stringent conditions on $r_i(t)$ that were used in other papers.
The restriction
of exponentially decay of $r_i(t)$ with $i$ (required in \cite{a1})
is replaced by much weaker condition
of convergence of the series $\sum_i g(i) r_i(t)$,
where the function $g(i)$ depends on
the characteristic which has to be calculated
using the truncated process (for example, the expected value
of the process). With respect to other intensities
 $\mu_{i}(t)$, $\lambda_{i}(t)$ and
$\beta_i(t)$
the only necessary assumption needed is that
any their linear combination has a single uniform upper bound.
In what follows we heavily rely on methodology developed in
\cite{gz04,z95,z08,z06}, which is based on the logarithmic norm
of linear operators and special transformations of the intensity
matrix governing the behaviour of the considered Markov process.

%In \cite{a1} it was not allowed for
%transition intensities to the state $0$
%to depend on the current state of the process.
%Here this restriction is not imposed.

The article is organized as follows. In the next
section the subclass of the birth and death
processes under consideration
and auxiliary results are introduced.
In section 3 it is explained how one
can obtain the ergodicity bounds.
Section 4 is devoted to method of truncation
that allows the
calculation of the limiting characteristics.
In section 5 one shows several examples of
how the obtained results can be
applied for the calculation of performance characteristics
of a specific queueing model.
Conclusion gives directions of future research.

\section{Description of the birth and death process}

Let the process $X\left(t\right)$, $t\geq 0$, be an inhomogeneous
continuous-time Markov chain with state space $\mathcal{X}=\{ 0, 1, 2 \dots \}$.
Transition, whenever it occurs from state $0$
can be to any state $i>0$, $i \in \mathcal{X}$. Transition from state
$i>0$ can be either to neighbouring state $\left(i-1\right)$ or $\left(i+1\right)$,
or to state $0$. All possible transition intensities are assumed to be non-random functions of time,
and may depend  (except for transition to the $0$)
on the state of the process.
Denote by $p_{ij}\left(s,t\right)=\Pr \left\{ X\left(t\right)=j\left| X\left(s\right)=i\right. \right\}$,
$i,j \ge 0, \;0\leq s\leq t$ transition probabilities of
$X\left(t\right)$ and by  $p_i\left(t\right)=\Pr \left\{ X\left(t\right) =i \right\}$
probability that Markov chain $X\left(t\right)$ is in state $i$ at time $t$.
Let ${\bf p}\left(t\right) = \left(p_0\left(t\right), p_1\left(t\right), \dots\right)^T$ be
probability distribution vector at time $t$.
Throughout the paper we assume that for $j \neq i$
\begin{equation}
\Pr\left(X\left( t+h\right) =j|X\left( t\right) =i\right)
=\begin{cases}\lambda_{i}\left( t\right)  h+\alpha_{ij}\left(t,
h\right), & \mbox { if }j =  i+1,\ i > 0, \cr \mu_{i}\left( t\right)
h+\alpha_{ij}\left(t, h\right), & \mbox { if }j =  i-1, \ i > 1, \cr
\beta_i\left( t\right)  h+\alpha_{ij}\left(t, h\right), & \mbox { if
}j = 0, \ i > 1, \cr r_j\left( t\right) h+\alpha_{ij}\left(t,
h\right), & \mbox { if } j \ge 1, \ i =0,  \cr \left(\mu_{1}\left(
t\right) + \beta_1\left( t\right)\right) h+\alpha_{ij}\left(t,
h\right), & \mbox { if }j =  i-1, \ i = 1, \cr \alpha_{ij}\left(t,
h\right) & \mbox { otherwise},
\end{cases}
\label{1001}
\end{equation}
\noindent where all  $\alpha_{i}(t,h)$ are $o(h)$ uniformly in $i$,
i. e., $\sup_i |\alpha_i(t,h)| = o(h)$.
Intensity functions $\beta_j\left( t\right)$ and $r_j\left( t\right)$, $j \ge 1$,
are henceforth called mass exodus
and mass arrivals intensities.

We assume that all intensity functions are
linear combinations of a finite number of locally integrable on
$[0,\infty )$ nonnegative functions. Then the corresponding intensity
matrix is
\begin{small}
$$
Q\left(t\right)=\left(
\begin{array}{cccccccc}
a_{00}\left(t\right) & r_1\left(t\right)  & r_2\left(t\right)   & r_3\left(t\right)  & r_4\left(t\right) & \ldots  & \ldots\\[3pt]
\beta_1\left(t\right) + \mu_1\left(t\right)   & a_{11}\left(t\right)  & \lambda_1\left(t\right)  & 0   & 0  & \ldots  & \ldots\\[3pt]
\beta_2\left(t\right)  & \mu_2\left(t\right)    & a_{22}\left(t\right)& \lambda_2\left(t\right)  & 0    &  \ldots &   \ldots\\[3pt]
\ldots&\ldots&\ldots&\ldots&\ldots&\ldots&\ldots\\[3pt]
\beta_j\left(t\right) & 0 &  \ldots & \mu_j\left(t\right)  &   a_{jj}\left(t\right) & \lambda_j\left(t\right)  & \ldots\\[3pt]
\vdots&\vdots&\vdots&\vdots&\vdots&\vdots&\ddots
\end{array}
\right),
$$
\end{small}

\noindent Let $a_{ij}\left(t\right) =  q_{ji}\left(t\right)$ for $j\neq i$ and
\begin{equation}
a_{ii}\left(t\right) =
-\sum_{j\neq i} a_{ji}\left(t\right) = -\sum_{j\neq i} q_{ij}\left(t\right).
\end{equation}
\noindent According to standard approach, which was applied in
\cite{gz04,z95,z08,z06}, we assume that the intensity matrix is
essentially bounded, i. e.
\begin{equation}
|a_{ii}\left(t\right)| \le L < \infty,
\label{0102-1}
\end{equation}
\noindent for almost all $t \ge 0$.

Probabilistic dynamics of
the considered process $X\left(t\right)$
is given by the
forward Kolmogorov system
\begin{equation}
\label{ur01}
\frac{d {\bf p} \left(t\right)}{dt}=A\left(t\right){\bf p} \left(t \right),
\end{equation}
\noindent where $A\left(t\right) = Q^T\left(t\right)$ is the transposed intensity matrix
of the process.
Throughout the paper by $\|\,\cdot\,\|$  we denote  the $l_1$-norm,
i. e.,  $\|{{\bf x} }\|=\sum_i |x_i|$, and $\|B\| = \sup_j \sum_i |b_{ij}|$
{for a matrix}  $B = \left(b_{ij}\right)_{i,j=0}^{\infty}$. Let $\Omega$ be the set all
stochastic vectors, i. e. $l_1$ -- vectors with non-negative coordinates
and unit norm.
Then we have $\|A\left(t\right)\| = 2\sup_{k}\left|a_{kk}\left(t\right)\right| \le 2 L $
for almost all $t \ge 0$. Hence, the operator function $A\left(t\right)$ from
$l_1$ into itself is bounded for almost all $t \ge 0$ and locally
integrable on $[0;\infty )$. Therefore, we can consider  (\ref{ur01})
as {a system of differential equations} in the space $l_1$ with bounded operator.

It is well known  (see \cite{DK}) that the Cauchy problem for
such system as  (\ref{ur01}) has a unique solution for
arbitrary initial condition and  ${\bf p} \left(s\right) \in \Omega$ implies
${\bf p} \left(t\right) \in \Omega$ for $t \ge s \ge 0$.

Denote by $E\left(t,k\right) = E\left\{X\left(t\right)\left|X\left(0\right)=k\right.\right\}$ the
expected value  (mean) of the process
$X\left(t\right)$ at moment $t$ under initial condition $X\left(0\right)=k$.

Recall that process $X\left(t\right)$ is called {\it weakly ergodic}, if
$\|{\bf p}^*\left(t\right)-{\bf p}^{**}\left(t\right)\| \to 0$ as $t \to \infty$ for any
initial conditions ${\bf p}^*\left(0\right), {\bf p}^{**}\left(0\right)$, where ${\bf
p}^*\left(t\right)$ and ${\bf p}^{**}\left(t\right)$ are the corresponding solutions of
(\ref{ur01}). Process $X\left(t\right)$ has  the {\it limiting mean}
$\varphi \left(t\right)$, if $ \lim_{t \to \infty }  \left(\varphi \left(t\right) -
E\left(t,k\right)\right) = 0$ for any $k$.

\section{Ergodicity bounds}

{In order to obtain new ergodicity bounds we apply the approach from
\cite{a1}. If for each instant $t$ one denotes by
$\beta_*\left(t\right)$ the greatest lower bound of
$\beta_i\left(t\right)$ with respect to $i$, i.e.}
\begin{equation}
\beta_*\left(t\right) = \inf_i \beta_i\left(t\right),\label{2001}
\end{equation}
\noindent then the forward Kolmogorov system (\ref{ur01})
can be rewritten in the following form:
\begin{equation}
\frac{d {\bf p} }{dt}=A^*\left( t\right) {{\bf p} }  +{\bf g} \left(t\right), \quad t\ge 0.
\label{eq112''}
\end{equation}
\noindent Here  ${\bf g}\left(t\right)=\left(\beta_*\left(t\right),0,0, \dots\right)^T$,
$A^*\left(t\right)=\left ( a_{ij}^*\left(t\right)\right)_{i,j=0}^{\infty}$, and
\begin{eqnarray}
a_{ij}^*\left(t\right) =   \left\{
\begin{array}{ccccccc}
a_{0j}\left(t\right) - \beta_*\left(t\right), & \mbox { if }  i= 0, \\
a_{ij}\left(t\right), & \mbox { otherwise}.
\end{array}
\right.  \label{1101}
\end{eqnarray}

{Further results heavily rely of the notion of the logarithmic norm of an operator
function.
Let $B\left( t\right) ,\ t\ge 0$ be a one-parameter family of
bounded linear operators on a Banach space ${\cal B}$ and let $I$
denote the identity operator. For a given $t\ge 0$, the number
\begin{eqnarray}
\gamma \left( B\left( t\right) \right)_{\cal B} =\lim\limits_{h\rightarrow
+0}\frac{%
\left\| I+hB\left( t\right) \right\| -1}h \label{lognorm}
\end{eqnarray}
is called the logarithmic norm of the operator $B\left( t\right)$\footnote{For further details on the logarithmic norm one can refer, for example, to \cite{gz04,dzp,z06}.}}.
If ${\cal B}$ is an $\left(N+1\right)-$dimensional vector space with $l_1$-
norm such that the operator $B\left(t\right)$ is given by the matrix
$B\left(t\right)=\left( b_{ij}\left(t\right)\right) _{i,j=0}^N$, $t\ge 0$, then the
logarithmic norm of $B\left(t\right)$ can be found explicitly:
$$
\gamma \left( B\left( t\right) \right) =\sup\limits_j\left( b_{jj}\left(
t\right) +\sum\limits_{i\neq j}\left| b_{ij}\left( t\right) \right| \right)
,\quad t\ge 0.
$$
On the other hand, the logarithmic norm of the operator $B\left(t\right)$ is
related to  the Cauchy operator $V\left(t,s\right)$ of the system
$$\frac{d{\bf x}}{dt}=B\left( t\right) {\bf x},\quad t\ge 0$$
in the following way:
$$
\gamma \left(B\left( t\right) \right)_{\cal B} =\lim\limits_{h\rightarrow
+0}\frac{%
\left\| V\left( t+h,t\right) \right\| -1}h, \quad t\ge 0.
$$
From the latter relation one can  deduce the  following bounds of the Cauchy
operator $V\left(t,s\right)$:
$$\left\| V\left( t,s\right) \right\|_{\cal B} \leq e^{\int\limits_s^t\gamma \left(
B\left( \tau \right) \right) \ d\tau },\quad 0\le s\le t.  $$ Here
we can find the exact value of the logarithmic norm of operator
function $A^*\left(t\right)$ , namely
\begin{equation}
\gamma \left(A^*\left(t\right)\right) = \sup_i \left(a_{ii}^*\left(t\right) + \sum_{j\neq
i}\left| a_{ji}^*\left(t\right)\right|\right) = -\beta_*\left(t\right). \label{cat03'}
\end{equation}

Let $U^*\left(t,s\right)$ be the Cauchy operator for equation (\ref{eq112''}).
Then $\|U^*\left(t,s\right)\| \le  e^{-\int_s^t \beta_*\left(\tau\right)\, d\tau}$,
and hence for any initial conditions ${{\bf p} }^{*}\left(0\right), {{\bf p} }^{**}\left(0\right)$
and any $t \ge 0$ we have
\begin{equation}
\left\|{{\bf p} }^{*}\left(t\right)-{{\bf p} }^{**}\left(t\right)\right\|\le e^{-\int\limits_0^t
\beta_*\left(\tau\right)\, d\tau} \left\|{{\bf p} }^{*}\left(0\right)-{{\bf p} }^{**}\left(0\right)\right\|,
\label{bound3'}
\end{equation}
\noindent Thus the following statement holds.

\begin{theorem} {\bf .} \label{t1}
Let the catastrophe rates be essential, i.e.,
\begin{equation}
\label{bet*}
\int_0^\infty \beta_*\left(t\right)\, dt = \infty.
\end{equation}
Then the process $X\left(t\right)$ is weakly ergodic (in the uniform operator
topology) and the following bound for the rate of convergence holds:
\begin{eqnarray}
\left\|{{\bf p}}^{*}\left(t\right)-{{\bf p}}^{**}\left(t\right)\right\|\le 2 e^{-\int\limits_0^t
\beta_*\left(\tau\right)\, d\tau}, \label{2011}\end{eqnarray} \noindent for any
initial conditions ${{\bf p}}^{*}\left(0\right), {{\bf p}}^{**}\left(0\right)$ and any $t \ge 0$.
\end{theorem}

Let now ${D}$ be such a diagonal matrix (${D} =
diag\left(d_0, d_1, d_2, \dots \right)$) such that the inequalities $1 = d_0 \le d_1 \le
\dots$ hold. Consider the corresponding space of sequences
$l_{1 {D}}=\left\{{\bf z} =\left(p_0,p_1,p_2,\ldots\right)\right\}$ such
that $\|{\bf z}\|_{1{D}}=\|{D} {\bf z}\|_1 <\infty$.
Then we can obtain the following estimate for the logarithmic norm
of the operator function $A^*\left(t\right)$ in the $l_{1{D}}$-norm:
\begin{eqnarray}
\gamma \left(A^*\left(t\right)\right)_{1{D}} = \gamma \left({D}
A^*\left(t\right) {D}^{-1}\right) = \sup_i \left(a_{ii}^*\left(t\right) +
\sum_{j\neq i}\left| \frac{d_j}{d_i}a_{ji}^*\left(t\right)\right|\right) = -
\beta_{**}\left(t\right), \label{0031'}
\end{eqnarray}
where $\beta_{**}\left(t\right) =  \inf_i \left(|a_{ii}^*\left(t\right)| - \sum_{j\neq
i}\left| \frac{d_j}{d_i}a_{ji}^*\left(t\right)\right|\right)$.

\begin{theorem}  {\bf .} \label{t2}
 Let there exist a sequence
$\{d_i\}$ such that
\begin{equation}
\int_0^\infty \beta_{**}\left(t\right)\, dt = \infty.
\end{equation}

Then  the following bound of the rate of convergence
holds:
\begin{eqnarray}
\left\|{{\bf p}}^{*}\left(t\right)-{{\bf p}}^{**}\left(t\right)\right\|_{1D}\le e^{-\int\limits_0^t
\beta_{**}\left(\tau\right)\, d\tau} \left\|{{\bf p}}^{*}\left(0\right)-{{\bf p}}^{**}\left(0\right)\right\|_{1D}, \label{201111}
\end{eqnarray}
\noindent for any initial conditions ${{\bf p}}^{*}\left(0\right), {{\bf p}}^{**}\left(0\right)$ and any $t \ge 0$.
\end{theorem}

Let $l_{1E}=\left\{{\bf z} =\left(p_0,p_1,p_2,\ldots\right)\right\}$ be a set
 of sequences such that $\|{\bf z}\|_{1E}=\sum_{k \ge 0} k
|p_k| <\infty$. Put $W = \inf_{k \ge 1}\frac{d_k}{k}$. Then $W\|{\bf
z}\|_{1E} \le \|{\bf z}\|_{1{D}}$, and we have the following
statement.

\begin{corollary}  {\bf .} \label{cor1}
 Let the conditions of Theorem  \ref{t2} hold and $W > 0$.
Then the  process $X\left(t\right)$ has the limiting mean, say
$\phi\left(t\right)=E\left(t,0\right)$, and the following bound of the rate of convergence
holds:
\begin{equation}
|E\left(t,j\right) - E\left(t,0\right)| \le  \frac{1+d_j}{W}e^{-\int\limits_0^t
\beta_{**}\left(\tau\right)\, d\tau}. \label{3012}
\end{equation}
\noindent for
any initial condition $j$ and any $t \ge 0$.
\end{corollary}

\section{Truncation bounds}

Consider the ``truncated'' process $X_{N}(t)$ on the state space
$E_{N} = \{0,1,\dots,N\}$ with the corresponding reduced intensity
matrix $A_N(t)$. Below we will identify the finite vector with
entries, say $ (a_0, a_1, \dots, a_{N})^{T}$ and the infinite vector
with the same first $N$ coordinates and others equal to zero.
Let us rewrite the system (\ref{eq112''}) as
\begin{equation}
\frac{d{\bf y}_N}{dt}=A^*(t){\bf y}_N + {\bf g}\left(t\right)+ \left(A^*_N\left(t\right)-A^*\left(t\right)\right) {\bf y}_N.
\end{equation}
and consider the corresponding ``truncated'' system
\begin{equation} \label{uss}
\frac{d{\bf y}_N}{dt}=A^*_N(t){\bf y}_N + {\bf g}_N\left(t\right).
\end{equation}
\noindent The solutions to the system (\ref{uss}) and (\ref{eq112''})  have the following form respectively
\begin{eqnarray}
{\bf p} \left(t\right) =  U^* \left(t, 0\right){\bf p} \left(0\right)  + \int\limits_0^t
U^* \left(t, \tau\right)   {\bf g}\left(\tau\right)  \, d\tau.
\label{01}
\\
 {\bf y}_N \left(t\right) = U^* \left(t, 0\right){\bf y}_N \left(0\right)  + \int\limits_0^t
U^* \left(t,
 \tau\right)   {\bf g}\left(\tau\right)  \, d\tau + \nonumber \\
 +  \int\limits_0^t
U^* \left(t,
 \tau\right)  \left(A^*_N\left(\tau\right)-A^*\left(\tau\right)\right) {\bf y}_N \left(\tau\right)  \, d\tau.
\label{011}
\end{eqnarray}

{Let the initial condition for both processes coincide, i.e.} ${\bf y}_N \left(0\right) = {\bf p} \left(0\right)$. Then in any norm one can write
\begin{equation}
\left\|{\bf p} \left(t\right) - {\bf y}_N \left(t\right)\right\| \le  \int\limits_0^t\left\|
U^* \left(t,
 \tau\right) \right\| \left\| \left(A^*_N\left(\tau\right)-A^*\left(\tau\right)\right) {\bf y}_N \left(\tau\right) \right\|\, d\tau.
\label{083}
\end{equation}

To evaluate the
Cauchy matrix  we use the logarithmic norm (\ref{lognorm}), and
we additionally assume that there are constants $M$ and $a$ such that
\begin{equation}
\left\| U^* \left(t, \tau\right) \right\| \le e^{-\int_\tau^t \beta^*\left(u\right) du} \le Me^{-a\left(t-\tau\right)}. \label{u*}
\end{equation}

We use the equality (\ref{201111}), and  also we will assume that there are constants $M_1$ and $a_1$ such that
\begin{eqnarray}
\|U^*\left(t, \tau\right)\|_{1D} \le e^{-\int_\tau^t \beta^{**} \left(u\right) du} \le M_1 e^{-a_1\left(t-\tau\right)}.
\end{eqnarray}

{Choose any constant $\theta$ such that $\theta \ge \sup_t {\beta^*\left(t\right)}$. Then one can bound
the solution of the system (\ref{eq112''}) in the $l_{1{D}}$-norm, introduced above, as follows}
\begin{eqnarray} \label{pD}
\|{\bf p} \left(t\right)\|_{1D}  = \!\!\!\!\!\!\!\!\!    && \| U^* \left(t, 0\right){\bf p} \left(0\right)  + \int\limits_0^t
U^* \left(t, \tau\right)   {\bf g}\left(\tau\right)  \, d\tau \| \le
\nonumber
\\ \le \!\!\!\!\!\!\!\!   &&
 \|U^*\left(t,0\right)\|_{1D}\|{\bf p} \left(0\right)\|_{1D} + \int_0^t \|U^*\left(t, \tau\right)\|_{1D}\| {\bf g} \left(\tau\right)\|_{1D} d\tau  \le
 \nonumber
 \\ \le \!\!\!\!\!\!\!\!   &&
M_1 e^{-a_1t}\|{\bf p} \left(0\right)\|_{1D} + \frac{M_1 \theta}{a_1} ,
\end{eqnarray}
\noindent Noticing that
$\|U^*_N\left(t,s\right)\|_{1D} \le \|U^*\left(t,0\right)\|_{1D}$ and  ${\bf g}_N \left(\tau\right) = {\bf g} \left(\tau\right)$
the solution of the system (\ref{uss})
in the $l_{1{D}}$-norm can be bounded as
\begin{eqnarray} \label{pND}
\|{\bf y}_N \left(t\right)\|_{1D} \le \!\!\!\!\! && \| U^*_N \left(t, 0\right){\bf y}_N \left(0\right) \| + \int\limits_0^t \|
U^*_N \left(t, \tau\right)   {\bf g}_N\left(\tau\right) \| \, d\tau  \le \nonumber
\\ \le && \!\!\!\!\!
 \|U^*\left(t,0\right)\|_{1D}\|{\bf p}\left(0\right)\|_{1D} + \int_0^t \|U^*\left(t, \tau\right)\|_{1D}\| {\bf g} \left(\tau\right)\|_{1D} d\tau  \le \nonumber
 \\ \le && \!\!\!\!\!
M_1 e^{-a_1t}\|{\bf p} \left(0\right)\|_{1D} + \frac{M_1 \theta}{a_1}.
\end{eqnarray}
Therefore
\begin{eqnarray} \label{pN}
p_N = \frac{d_N p_N}{d_N}  \le \frac{\|{\bf y}_N \left(t\right)\|_{1D}}{d_N} \le \frac{M_1 \theta}{a_1 d_N}  + \frac{M_1 e^{-a_1t}}{d_N} \|{\bf p} \left(0\right)\|_{1D}.
\end{eqnarray}
Next, noticing that the difference $A^*\left(\tau\right)- A^*_N\left(\tau\right)$
is equal to
\begin{small}
\begin{eqnarray}
A^*\left(\tau\right)- A^*_N\left(\tau\right)=\left(
\begin{array}{ccccccccc}
-\sum_{n \ge N+1} r_n & 0  & 0   & 0  & 0 & \ldots  & \beta_{N+1}-\beta^* & \ldots & \ldots \\
0 & 0  & 0  & 0   & 0  & \ldots  & 0 & \ldots & \ldots\\
\ldots&\ldots&\ldots&\ldots&\ldots&\ldots&\ldots&\ldots& \ldots \\
0  & 0 & \ldots & \ldots  & 0& -\lambda_N  & \mu_{N+1}    &  0 &   \ldots \\
r_{N+1} & 0 & 0 & \ldots & \ldots & \lambda_N& a_{N+1,N+1}& \mu_{N+2} & 0 \\
\ldots&\ldots&\ldots&\ldots&\ldots&\ldots&\lambda_{N+1}&\ldots &\ldots
\end{array}
\right),
\end{eqnarray}
\end{small}
and for the $\left({\bf A}^*\left(\tau\right)-{\bf A}^*_N\left(\tau\right)\right){\bf y}_N $
we have
\begin{eqnarray}
\left({\bf A}^*\left(\tau\right)-{\bf A}^*_N\left(\tau\right)\right){\bf y}_N  = \left(
\begin{array}{c}
- \sum_{n \ge N+1} r_n y_0 \\
0 \\
\ldots\\
-\lambda_N y_N \\
r_{N+1} y_0 + \lambda_N y_N \\
r_{N+2} y_0 \\
\ldots
\end{array}
\right).
\end{eqnarray}
then, using bound (\ref{pN}), we have the following bound:
\begin{eqnarray}\label{oan}
\left\| \left( A^*_N\left(\tau\right)- A^*\left(\tau\right)\right)
{\bf y}_N \left(\tau\right) \right\| \le \!\!\!\!\! && 2\sum_{n \ge N+1} r_n p_0 +2 \lambda_N p_N  \le
\\ \nonumber \le \!\!\!\!\! &&
 R_{N+1}+ \frac{2 M_1 L}{d_N}\left( \frac{ \theta}{a_1}  +  e^{-a_1t}\|{\bf p} \left(0\right)\|_{1D} \right) ,
\end{eqnarray}
\noindent where $R_{N+1}$  is the remainder of the convergent series $\sum_{n}^\infty r_{n}$.

Therefore, under assumptions (\ref{083}), (\ref{u*}),
one has the following bound of truncation from (\ref{oan}):
\begin{eqnarray*}
\left\|{\bf p} \left(t\right) - {\bf y}_N \left(t\right)\right\| \le \quad
\quad \quad \quad \quad \quad \quad \quad \quad \quad \quad \quad \quad
\quad \quad \quad \quad \quad \quad \quad \quad \quad \quad \quad \quad \\
\le  \left( MR_{N+1}+ \frac{2M M_1 L}{d_N}\left( \frac{ \theta}{a_1}  +  e^{-a_1t}\|{\bf p} \left(0\right)\|_{1D} \right)\right) \int\limits_0^t
e^{-a\left(t-\tau\right)}\, d\tau  \le \\ \le
\frac{MR_{N+1}}{a}+ \frac{2M M_1 L}{ad_N}\left( \frac{\theta}{a_1}  +  e^{-a_1t}\|{\bf p} \left(0\right)\|_{1D}\right).
\end{eqnarray*}

\begin{theorem} {\bf .} \label{t3}
Let the conditions of Theorems \ref{t1} and \ref{t2} hold.
Then the following bound   of truncation holds for any
initial conditions ${{\bf p}}\left(0\right)={{\bf y}}_N\left(0\right)$ and any $t \ge 0$
\begin{eqnarray}
\left\|{\bf p} \left(t\right) - {\bf y}_N \left(t\right)\right\| \le
\frac{MR_{N+1}}{a}+ \frac{2M M_1 L}{ad_N}\left( \frac{ \theta}{a_1}  +  e^{-a_1t}\|{\bf p} \left(0\right)\|_{1D}\right).
\end{eqnarray}
\end{theorem}

In order to find the bound for the mean of the truncated process
consider the coordinates of the vector ${{\bf p}}\left( t \right)$.
We have
\begin{eqnarray} \label{ddd}
p_{N+1}+ \!\!\! && \!\!\!\!\!\! 2p_{N+2}+3p_{N+3}+\cdots \le \nonumber \\
\le && \!\!\!\!\!\! \frac{1}{d_{N+1}}d_{N+1}p_{N+1}+\frac{2}{d_{N+2}}d_{N+2}p_{N+2}+\frac{3}{d_{N+3}}d_{N+3}p_{N+3}+\cdots \le  \nonumber  \\ \le && \!\!\!\!\!\!
\left\|{\bf p} \left(t\right) \right\| _{1D} \sup_{i \ge 1}{\frac{i}{d_{N+i}}}.
\end{eqnarray}
{From \eqref{ddd} we get the bound for the difference in means of
the original and the truncated process:}
\begin{eqnarray*}
\left| E(t,k)-E_N(t,k) \right| \le
\left\|{\bf p} \left(t\right) - {\bf y}_N \left(t\right)\right\|_{1E} < \\
 < \left| p_0 - p_{N0} \right| +  \left| p_1-p_{N1} \right| + 2\left| p_2-p_{N2} \right| + \cdots +\\+N\left| p_N-p_{NN} \right| +\left( N+1 \right) p_{N+1}+\left(N+2\right)p_{N+2}+\left(N+3\right)p_{N+3}+\cdots \le \\ \le N\left\|{\bf p} \left(t\right) - {\bf y}_N \left(t\right)\right\|+ p_{N+1}+2p_{N+2}+3p_{N+3}+\cdots
\end{eqnarray*}
From relations (\ref{pD}) and  (\ref{ddd}) one can obtain the following corollary.

\begin{corollary}{\bf .} \label{cor2}
Let the conditions of Theorem \ref{t3} hold and assume the series $\sum_{n} n r_n$ converges.
Then for any
initial conditions $X(0)=X_N(0)=k$ and any $t \ge 0$
the following truncation bound holds:
\begin{eqnarray}
\left| E(t,k)-E_N(t,k) \right| \le \frac{N MR_{N+1}}{a}+ \quad\quad\quad\quad\quad\quad\quad\quad\quad\quad \quad\quad\quad\quad\quad\quad\quad\quad \\  +  \left(\frac{2N M M_1 L}{ad_N} + \sup_{i \ge 1}\frac{i}{d_{N+i}} \right)\left( \frac{ \theta}{a_1}  +  e^{-a_1t}\|{\bf p} \left(0\right)\|_{1D}\right),\nonumber
\end{eqnarray}
where $\|{\bf p}(0)\|_{1D}=d_k$.
\end{corollary}

\section{Examples}

Efficiency of the bounds obtained for the process
$X(t)$ in the previous sections will illustrated in the queueing theory context.
Specifically we consider
$M_t/M_t/S$ queue with catastrophes and bulk arrivals when
empty when intensities are periodic functions of time
which can be described by process $X(t)$.
In each example it is shown how to find approximations for the limiting
value of the mean number of customers in the system and limiting value of empty system
with given error. For convenience we first give detailed description of the system
and then proceed to examples.
Queueing system consists of single infinite capacity queue
and $S$ servers. Two flows of customers arrive at
the system: flow of ordinary customers and flow of catastrophes.
If at time $t$ there is at least one customer in the system
then new arrivals of ordinary customers happen according
to inhomogeneous Poisson process with intensity $\lambda(t)$.
But if at time $t$ the system is empty ordinary customers arrive
in bulk (or groups) in accordance with a inhomogeneous Poisson process
of intensity $r(t)$. The size of arriving group is a random variable
with probability distribution $g_n$, $n=1,2, \dots$, having finite mean.
The sizes and interarrival times of successive arriving groups are
stochastically independent. Let $r_n(t)=g_n r(t)$. Each ordinary customer upon arrival
occupies one place in the queue and waits for service. Whenever server
becomes free customer from the queue (if there is any) enters server and
get served according to exponential distribution
with intensity $\mu(t)$ (service discipline is unimportant
and for certainty one can consider that customer are served in FIFO manner).
Additional inhomogeneous Poisson flow of catastrophes of intensity $\beta_n(t)$
arrives at the system. If arriving customer of this flow
finds the system busy it removes all customers from the system
and leaves it. Otherwise it has no effect on it.

In order to illustrate the behaviour of
system's performance characteristics
we will consider several special cases: when
there are few servers in the system and when there
are many servers in the system.
We will consider three performance characteristics:
the limiting probability $p_0(t)$ of the empty system, the
limiting probability $\Pr(X(t)  \le S )=\sum_{i=0}^S p_i(t)$ of the empty queue,
and the limiting mean $E(t,0)$ if initially the system was empty.

\subsection{Case when S=3}

Let the number of servers in the system be equal to $S=3$. We will
consider two examples which differ in the values for the intensities
$r_n(t)$.

\noindent { \bf Example 1.}
%\begin{example}
Let the intensities have the form
\begin{eqnarray*}
\lambda_n\left(t\right)&=& \lambda\left(t\right) = 1+\sin 2\pi t,
\\
\mu_n\left(t\right)&=&
\min\left(n,S\right)\mu\left(t\right)=\min\left(n,3\right)\left(1-\sin 2\pi t\right),
\\
\beta_n\left(t\right)&=& 2+\cos2\pi t+ \frac{1}{n},
\\
r_n\left(t\right)&=&\frac{1+\sin2\pi t}{4^n}.
\end{eqnarray*}
%\end{example}

This is almost the same example as example 1.1 in \cite{a1},
except for the fact that here the intensities $\beta_n(t)$ are dependent on $n$,
but in \cite{a1} such dependency was not allowed.
Notice also that here the values of $r_n(t)$ decay exponentially.
%
%But with the growth of $n$ the value of
%$\beta_n(t)$ becomes independent of $n$.
Due to the fact that \eqref{bet*} holds, then according to
\textit{Theorem} \ref{t1} the process $X\left(t\right)$ is
exponentially weakly ergodic. Hence $X\left(t\right)$ has the
limiting $1-$periodic regime and the correspondent limiting
$1-$periodic mean, see details in \cite{z06,a3}. The following bound
for the rate of convergence holds:
\begin{eqnarray}
\left\|{\bf p}^{*}\left(t\right)-{\bf
p}^{**}\left(t\right)\right\|\le 2 e^{-\int\limits_0^t
\beta\left(\tau\right)\, d\tau} \le 4 e^{-2t} ,
\label{71}\end{eqnarray} \noindent for any initial conditions ${\bf
p}^{*}\left(0\right), {\bf p}^{**}\left(0\right)$ and any $t \ge 0$.

Let $d_n=2^{n}$. One has that $L=12$, $M = e^{1/\pi}$, $a=2$,
$\beta_{**}\left( t \right) = \frac{3}{2} + \cos2\pi t -
\frac{3}{2}\sin2\pi t$, $a_1=1.5$ and  $M_1 = e^{2.5/ \pi}$, and
Corollary 1 implies the following bound
\begin{equation}
|E\left(t,j\right) - E\left(t,0\right)| \le  3\left(d_j+1\right)
e^{-1.5t} \label{72}
\end{equation}
\noindent for any initial condition $j$ and any $t \ge 0$.

Then using \textit{Theorem} \ref{t3} and \textit{Corollary}
\ref{cor2} with truncation error $10^{-5}$ for $N=30$ and $t \in
[10,11]$ one obtains the following bounds of truncation:
\begin{eqnarray}
\left\|{\bf p} \left(t\right) - {\bf y}_N \left(t\right)\right\| \le 10^{-7},
\end{eqnarray}
\begin{eqnarray}
\left|E\left(t,0\right)- E_N\left(t,0\right)\right| \le  3\cdot
10^{-6}.
\end{eqnarray}

 In Fig.1--3 one can see the behaviour
of the limiting probabilities $p_0(t)$ and $\Pr(X(t)  \le S )$,
and the limiting mean $E(t,0)$ for different values of $t$.

\bigskip

\noindent { \bf Example 2.}
Let the intensities have the form
\begin{eqnarray*}
\lambda_n\left(t\right)&=& \lambda\left(t\right) = 1+\sin 2\pi t,
\\
\mu_n\left(t\right)&=&
\min\left(n,S\right)\mu\left(t\right)=\min\left(n,3\right)\left(1-\sin 2\pi t\right),
\\
\beta_n\left(t\right)&=& 2+\cos2\pi t+ \frac{1}{n},
\\
r_n\left(t\right)&=&\frac{1+\sin2\pi t}{n^{10}}.
\end{eqnarray*}

Notice that though the series of $\sum_n r_n(t)$ though converges rapidly its $n$-th term does not decay exponentially because (from starting from a certain value) $r_n(t)$ is always greater than the corresponding term $r_n(t)$ in the previous example.
Due to the fact that \eqref{bet*} holds, then according to
\textit{Theorem} \ref{t1} the process $X\left(t\right)$ is
exponentially weakly ergodic. Hence $X\left(t\right)$ has the
limiting $1-$periodic regime and the correspondent limiting
$1-$periodic mean. The following bound
for the rate of convergence holds:
\begin{eqnarray}
\left\|{\bf p}^{*}\left(t\right)-{\bf p}^{**}\left(t\right)\right\|\le 2 e^{-\int\limits_0^t
\beta\left(\tau\right)\, d\tau} \le 4 e^{-2t} , \label{71}
\end{eqnarray}
\noindent for any initial conditions ${\bf p}^{*}\left(0\right), {\bf p}^{**}\left(0\right)$
and any $t \ge 0$.

Let $d_n=\frac{3}{2}^{n}$, if $n < 100$, and
$d_n=\frac{3}{2}^{100}\frac{n+1}{100}$, if $n \ge 100$. One has that
$L=12$, $M = e^{1/\pi}$, $a=2$, $\beta_{**}\left( t \right) =
\frac{5}{6} + \cos2\pi t - \frac{5}{6}\sin2\pi t$, $a_1=0.8$ and
$M_1 = e^{1.8/ \pi}$, and Corollary 1 implies the following bound
\begin{equation}
|E\left(t,j\right) - E\left(t,0\right)| \le  2\left(d_j+1\right)
e^{-0.8t} \label{72}
\end{equation}
\noindent for any initial condition $j$ and any $t \ge 0$.

Then using \textit{Theorem} \ref{t3} and \textit{Corollary}
\ref{cor2} with truncation error $10^{-5}$ for $N=55$ and $t \in [17,18]$ one obtains the following bounds of truncation:
\begin{eqnarray}
\left\|{\bf p} \left(t\right) - {\bf y}_N \left(t\right)\right\| \le 3\cdot 10^{-8},
\end{eqnarray}
\begin{eqnarray}
\left|E(t,0) - E_N(t,0)\right| \le  2\cdot 10^{-6}.
\end{eqnarray}

The behaviour
of the limiting probabilities $p_0(t)$ and $\Pr(X(t)  \le S )$,
and the limiting mean $E(t,0)$ for different values of $t$
one can see in in Fig.4--6.

\bigskip

\subsection{Case when S=20}

The following two examples show
what changes in the behaviour of system's performance characteristics
with the change in the number of servers, arrival intensities.

\noindent { \bf Example 3.}
Let the intensities have the form
\begin{eqnarray*}
\lambda_n\left(t\right)&=& \lambda\left(t\right) = 15\left( 1+\sin 2\pi t \right),
\\
\mu_n\left(t\right)&=&
\min\left(n,S\right)\mu\left(t\right)=\min\left(n,20\right)\left(1-\sin 2\pi t\right),
\\
\beta_n\left(t\right)&=& 2+\cos2\pi t+ \frac{1}{n},
\\
r_n\left(t\right)&=&\frac{1+\sin2\pi t}{4^n}.
\end{eqnarray*}

Due to the fact that \eqref{bet*} holds, then according to
\textit{Theorem} \ref{t1} the process $X\left(t\right)$ is
exponentially weakly ergodic. Hence $X\left(t\right)$ has the
limiting $1-$periodic regime and the correspondent limiting
$1-$periodic mean. The following bound
for the rate of convergence holds:
\begin{eqnarray}
\left\|{\bf p}^{*}\left(t\right)-{\bf p}^{**}\left(t\right)\right\|\le 2 e^{-\int\limits_0^t
\beta\left(\tau\right)\, d\tau} \le 4 e^{-2t} ,
\end{eqnarray}
\noindent for any initial conditions ${\bf p}^{*}\left(0\right), {\bf p}^{**}\left(0\right)$
and any $t \ge 0$.

Let $d_n=\left(1+\frac{1}{8}\right)^{n}$. One has that $L=74$, $M =
e^{1/\pi}$, $a=2$, $\beta_{**}\left( t \right) = \frac{559}{2520} +
\cos2\pi t - \frac{143}{72}\sin2\pi t$, $a_1=0.2$ and  $M_1 = e^{3/
\pi}$, and Corollary 1 implies the bound
\begin{equation}
|E\left(t,j\right) - E\left(t,0\right)| \le  3\left(d_j+1\right)
e^{-0.2t} \label{72}
\end{equation}
\noindent for any initial condition $j$ and any $t \ge 0$.

Then for exponentially decaying values of $r_n(t)$
using \textit{Theorem} \ref{t3} and corollary \ref{cor2} with
truncation error $10^{-5}$ for $N=220$ and $t \in [60,61]$
one obtains the following bounds of truncation:
\begin{eqnarray}
\left\|{\bf p} \left(t\right) - {\bf y}_N \left(t\right)\right\| \le 3\cdot 10^{-8},
\end{eqnarray}
\begin{eqnarray}
\left|E(t,0) - E_N(t,0)\right| \le  6\cdot 10^{-6}.
\end{eqnarray}

Fig. 7--9 show the behaviour of $p_0(t)$, $\Pr(X(t)  \le S )$
and $E(t,0)$ for different values of $t$.

\bigskip

\noindent { \bf Example 4.}
Let the intensities have the form
\begin{eqnarray*}
\lambda_n\left(t\right)&=& \lambda\left(t\right) = 15\left(1+\sin 2\pi t\right),
\\
\mu_n\left(t\right)&=&
\min\left(n,S\right)\mu\left(t\right)=\min\left(n,20\right)\left(1-\sin 2\pi t\right),
\\
\beta_n\left(t\right)&=& 2+\cos2\pi t+ \frac{1}{n},
\\
r_n\left(t\right)&=&\frac{1+\sin2\pi t}{n^{10}}.
\end{eqnarray*}

Due to the fact that \eqref{bet*} holds, then according to
\textit{Theorem} \ref{t1} the process $X\left(t\right)$ is
exponentially weakly ergodic. Hence $X\left(t\right)$ has the
limiting $1-$periodic regime and the correspondent limiting
$1-$periodic mean. The following bound
for the rate of convergence holds:
\begin{eqnarray}
\left\|{\bf p}^{*}\left(t\right)-{\bf p}^{**}\left(t\right)\right\|\le 2 e^{-\int\limits_0^t
\beta\left(\tau\right)\, d\tau} \le 4 e^{-2t} , \label{71}
\end{eqnarray}
\noindent for any initial conditions ${\bf p}^{*}\left(0\right), {\bf p}^{**}\left(0\right)$
and any $t \ge 0$.

Let $d_n=\frac{9}{8}^{n}$, if $n < 200$, and
$d_n=\frac{9}{8}^{200}\frac{n+1}{200}$, if $n \ge 200$. One has that
$L=74$, $M = e^{1/\pi}$, $a=2$, $\beta_{**}\left( t \right) =
\frac{17}{72} + \cos2\pi t - \frac{143}{72}\sin2\pi t$, $a_1=0.2$
and  $M_1 = e^{3/ \pi}$, and Corollary 1 implies the bound
\begin{equation}
|E\left(t,j\right) - E\left(t,0\right)| \le  3\left(d_j+1\right)
e^{-0.2t} \label{72}
\end{equation}
\noindent for any initial condition $j$ and any $t \ge 0$.

 Then using Theorem \ref{t3} and corollary \ref{cor2} with truncation error $10^{-5}$ for $N=220$ and $t \in [56,57]$ one obtains the following bounds of truncation:
\begin{eqnarray}
\left\|{\bf p} \left(t\right) - {\bf y}_N \left(t\right)\right\| \le 3\cdot 10^{-8},
\end{eqnarray}
\begin{eqnarray}
\left|E(t,0) - E_N(t,0)\right| \le  6\cdot 10^{-6}.
\end{eqnarray}

The behaviour of $p_0(t)$, $\Pr(X(t)  \le S )$
and $E(t,0)$ for different values of $t$
is shown in Fig. 10--12.

\section{Conclusion}

In this paper we have obtained
ergodicity and truncation bounds
for a class of inhomogeneous birth and
death process with an additional arrival from/to origin
under relaxed conditions on the transitions from the origin.
Previously these bounds were known only for
exponentially decaying transitions
from the origin.
The obtained results are especially accurate in cases when
the arrival rate is not very high and become worse
with its increase. The study of the boundaries for
the arrival rate which lead to accurate estimates
as well as the development of new approaches to
deal with high load values is a promising direction of research.

\bigskip

\section*{Acknowledgements}

This work was supported by Russian Scientific
Foundation   (Grant No. 14-11-00397).

\begin{figure}[!b]
 \label{fig1n}
\begin{center}
 \includegraphics[width=11cm]{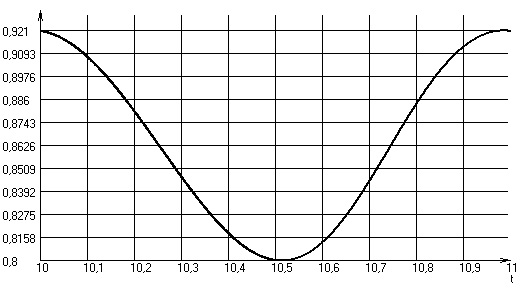}
 \caption{Case $S=3$. $r_n(t)$ decay exponentially. Approximation of the limiting probability of
empty queue $p_0(t)$ on $[10,11]$.}
\end{center}
\end{figure}

\begin{figure}[!b]
\begin{center}
 \includegraphics[width=11cm]{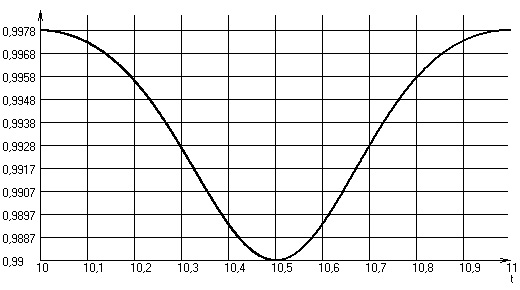}
 \caption{Example 1. Approximation of the limiting probability  $\Pr(X(t)  \le S )$ on $[10,11]$.}
\end{center}
 \label{fig2n}
\end{figure}

\begin{figure}[!b]
\begin{center}
 \includegraphics[width=11cm] {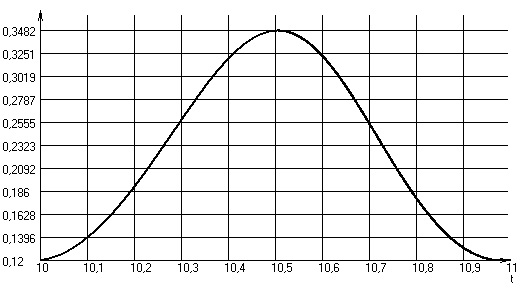}
 \caption{Example 1. Approximation of the limiting mean $E(t,0)$ on $[10,11]$.}
\end{center}
 \label{fig3n}
\end{figure}

\begin{figure}[!b]
\begin{center}
 \includegraphics[width=11cm]{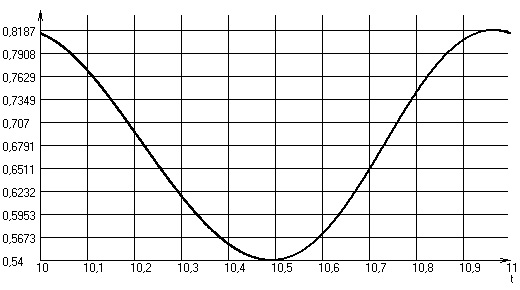}
 \caption{Example 2. Approximation of the limiting probability of
empty queue $Pr\left(X\left(t\right)=0\right)$ on $[17,18]$.}
\end{center}
 \label{fig4}
\end{figure}

\begin{figure}[!b]
\begin{center}
 \includegraphics[width=11cm]{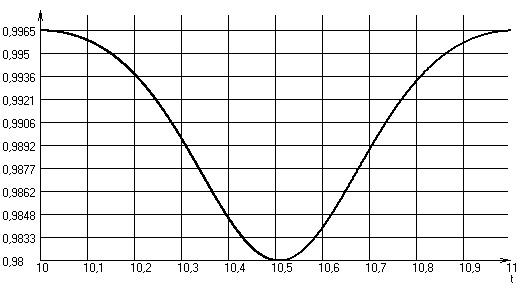}
 \caption{Example 2. Approximation of the limiting probability $\Pr\left(X\left(t\right)  \le S \right)$ on $[17,18]$.}
\end{center}
 \label{fig5}
\end{figure}

\begin{figure}[!b]
\begin{center}
 \includegraphics[width=11cm] {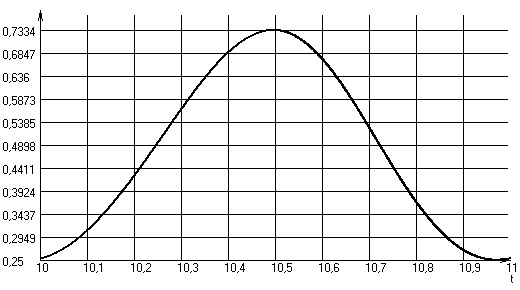}
 \caption{Example 2. Approximation of the limiting mean $E(t,0)$ on $[17,18]$.}
\end{center}
 \label{fig6}
\end{figure}

\begin{figure}[!b]
\begin{center}
\includegraphics[width=11cm]{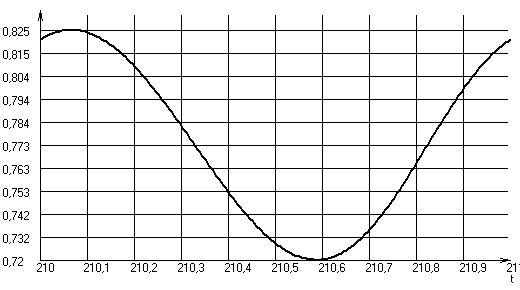}%
 \caption{Example 3. Approximation of the limiting probability of
empty queue $Pr\left(X\left(t\right)=0\right)$ on $[60,61]$.}
\end{center}
 \label{fig7}
\end{figure}

\begin{figure}[!b]
\begin{center}
\includegraphics[width=11cm]{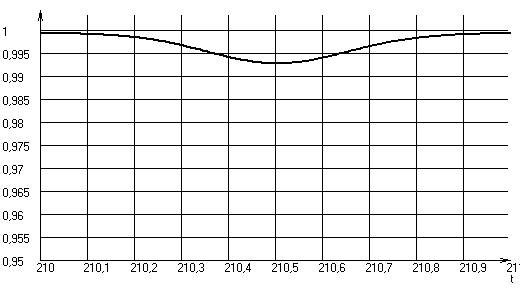}
 \caption{Example 3. Approximation of the limiting probability $\Pr\left(X\left(t\right)  \le S \right)$ on $[60,61]$.}
\end{center}
 \label{fig8}
\end{figure}

\begin{figure}[!b]
\begin{center}
 \includegraphics[width=11cm] {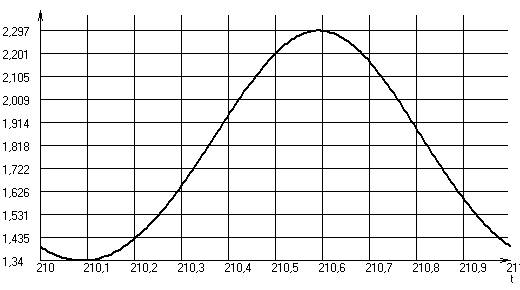}
 \caption{Example 3. Approximation of the limiting mean $E(t,0)$ on $[60,61]$.}
\end{center}
 \label{fig9}
\end{figure}

\begin{figure}[!b]
\begin{center}
 \includegraphics[width=11cm]{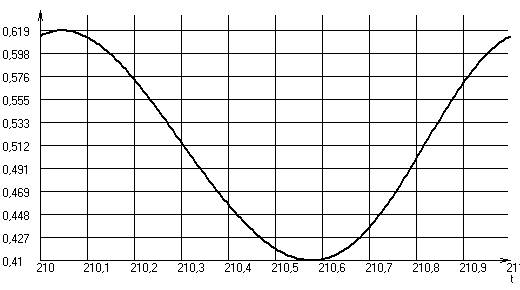}
 \caption{Example 4. Approximation of the limiting probability of
empty queue $\Pr\left(X\left(t\right)=0\right)$ on $[56,57]$.}
\end{center}
 \label{fig10}
\end{figure}

\begin{figure}[!b]
\begin{center}
 \includegraphics[width=11cm]{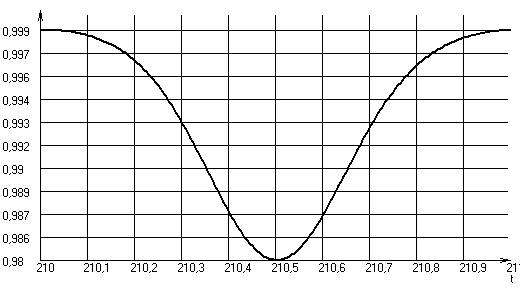}
 \caption{Example 4. Approximation of the limiting probability $\Pr\left(X\left(t\right)  \le S \right)$ on $[56,57]$.}
\end{center}
 \label{fig11}
\end{figure}

\begin{figure}[!b]
\begin{center}
\includegraphics[width=11cm] {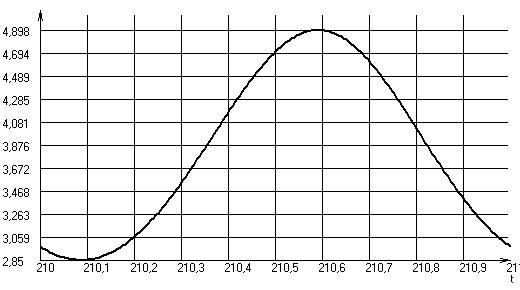}
 \caption{Example 4. Approximation of the limiting mean $E(t,0)$ on $[56,57]$.}
\end{center}
 \label{fig12}
\end{figure}

\end{document}